\documentstyle[12pt, epsf, epsfig, amssymb, amstext, amstex,amsthm, amsfonts]{article}
\input xy
\input amssymb.sty
\xyoption{all}
\begin{document}

\newtheorem{thm}{Theorem}[section]
\newtheorem{prop}[thm]{Proposition}
\newtheorem{cor}[thm]{Corollary}
\newtheorem{lem}[thm]{Lemma}
\newtheorem{conj}[thm]{Conjecture}
\newtheorem{exa}[thm]{Example}
\newtheorem{defn}[thm]{Definition}
\newtheorem{clm}[thm]{Claim}
\newtheorem{eex}[thm]{Exercise}
\newtheorem{obs}[thm]{Observation}
\newtheorem{note}[thm]{Notation}
 
\newcommand{\ben}{\begin{enumerate}}
\newcommand{\een}{\end{enumerate}}
\newcommand{\blem}{\begin{lem}}
\newcommand{\elem}{\end{lem}}
\newcommand{\bcl}{\begin{clm}}
\newcommand{\ecl}{\end{clm}}
\newcommand{\bthm}{\begin{thm}}
\newcommand{\ethm}{\end{thm}}
\newcommand{\bpr}{\begin{prop}}
\newcommand{\epr}{\end{prop}}
\newcommand{\bco}{\begin{cor}}
\newcommand{\eco}{\end{cor}}
\newcommand{\bcon}{\begin{conj}}
\newcommand{\econ}{\end{conj}}
\newcommand{\bde}{\begin{defn}}
\newcommand{\ede}{\end{defn}}
\newcommand{\bex}{\begin{exa}}
\newcommand{\eexa}{\end{exa}}
\newcommand{\bexe}{\begin{exe}}
\newcommand{\eexe}{\end{exe}}
\newcommand{\bobs}{\begin{obs}}
\newcommand{\eobs}{\end{obs}}
\newcommand{\bnote}{\begin{note}}
\newcommand{\enote}{\end{note}}

\newcommand{\fg}{\Pi _1(D-K,u)}
\newcommand{\Z}{{\Bbb Z}}
\newcommand{\C}{{\Bbb C}}
\newcommand{\R}{{\Bbb R}}
\newcommand{\Q}{{\Bbb Q}}
\newcommand{\F}{{\Bbb F}}
\newcommand{\N}{{\Bbb N}}

\newcommand{\fnref}[1]{~(\ref{#1})}
\newenvironment{emphit}{\begin{itemize} \em}{\end{itemize}}
\begin{center}
\Large{\bf {Graph Theoretic Method for Determining non- Hurwitz Equivalence in the Braid Group and Symmetric group}}\\ 
\vspace{7mm}
\large{T. Ben-Itzhak and M. Teicher}
\footnote{Partially supported by the Emmy Noether Research Institute for Mathematics
and the Minerva Foundation of Germany and to the Excellency Center "Group
Theoretic Methods in the Study of Algebraic Varieties"  of the Israel Science Foundation and by EAGER (European network in Algebraic Geometry).\\}
\end{center}


\vspace{15mm}
ABSTRACT. Motivated by the problem of Hurwitz equivalence of $\Delta ^2$ factorization in the braid group, we address the problem of Hurwitz equivalence in the symmetric group, obtained by projecting the $\Delta ^2$ factorizations into $S_n$. We get $1_{S_n}$ factorizations with transposition factors. Looking at the transpositions as the edges in a graph, we show that two factorizations are Hurwitz equivalent if and only if their graphs have the same weighted connected components. 
The main result of this paper will help us to compute the BMT invariant presented in \cite{KuTe} or \cite{Te}. The graph structure gives a weaker but very easy to compute invariant to distinguish between diffeomorphic surfaces which are not deformation of each other.

\section{Definitions}


\bde \label{BraidGroup}
Braid Group $B_n$
\ede
\noindent $B_n$ is the group generated by $\sigma _1,...,\sigma _{n-1}$ with the following relations:
$$\sigma _i \sigma _j = \sigma _j \sigma _i \quad |i-j|>1$$
$$\sigma _i \sigma _{i+1} \sigma _i = \sigma _{i+1} \sigma _i \sigma _{i+1}.$$


\bde \label{HT}
Half Twist
\ede
\noindent Let $H \in B_n$, we say that $H$ is a half twist if $H = P \sigma _i P ^{-1}$ for some $ 1 \leq i \leq n-2$ and $P \in B_n$.

\bde
Hurwitz move on $G^m$ ($R_k, R_k ^{-1}$)
\ede
\noindent Let $G$ be a group, $\overrightarrow{t}=(t_1,...,t_m) \in G^m$. We say that $\overrightarrow{s}=(s_1,...,s_m) \in G^m$ is obtained from $\overrightarrow{t}$ by the Hurwitz move $R_k$ (or $\overrightarrow{t}$ is obtained from $\overrightarrow{s}$ by the Hurwitz move $R_k ^{-1}$) if 
$$ s_i = t_i \quad \text{ for } i \not= k, k+1,$$   
$$ s_k = t_k t_{k+1} t_k ^{-1}, \quad s_{k+1} = t_k.$$
 

\bde
Hurwitz move on factorization 
\ede
\noindent Let $G$ be a group and $t \in G$. Let $t = t_1 \cdots t_m = s_1 \cdots s_m$ be two factorized expressions of $t$. We say that $s_1 \cdots s_m$ is obtained from $t_1 \cdots t_m$ by the Hurwitz move $R_k$ if $(s_1,...,s_m)$ is obtained from $(t_1,...,t_m)$ by the Hurwitz move $R_k$.\\


\bde
Hurwitz equivalence of factorization
\ede
\noindent The factorizations $s_1 \cdots s_m$, $t_1 \cdots t_m$\label{formin2} are Hurwitz equivalent if they are obtained from each other by a finite sequence of Hurwitz moves. The notation is  $t_1 \cdots t_m \overset{HE}{\backsim} s_1 \cdots s_m$.\label{formin1} \\


\section{Projecting to $S_n$}

Let $\phi : B_n \rightarrow S_n$ be the natural homomorphism to $S_n$, given by $\phi(b) \rightarrow \pi _b$ where $\pi _b$ 
is the permutation given by the strings of $b$. In terms of definition \ref{BraidGroup}, $\phi : B_n \rightarrow S_n$ is defined by $\phi (\sigma _i) = (i,i+1) \quad 1 \leq i \leq n-1$.


\bpr \label{Project}
Let $b_1 \cdots b_m$, $r_1 \cdots r_m$ be two factorizations in $B_n$ s.t. 
$b_1 \cdots b_m \overset{HE}{\backsim} r_1 \cdots r_m$, then, 
$\phi (b_1) \cdots \phi (b_m) \overset{HE}{\backsim} \phi (r_1) \cdots \phi (r_m)$.
\epr

\noindent
{\bf Proof:} It is sufficient to show for the case where the factorization $r_1 \cdots r_m$ is obtained from 
$b_1 \cdots b_m$ by a single Hurwitz move, $R_i$, and therefore, $r_i = b_i b_{i+1} {b_i}^{-1}$ and $r_{i+1} = b_i$
and $r_k = b_k$ if $k \neq i,i+1$.\\
By performing $R_i$ on $\phi (b_1) \cdots \phi (b_i) \cdot \phi (b_{i+1}) \cdots \phi (b_m)$ we get 
$$\phi (b_1) \cdots \phi (b_i) \phi (b_{i+1}) \phi ({b_i})^{-1} \cdot \phi (b_i) \cdots \phi (b_m)$$
which is equal to,
$$\phi (b_1) \cdots \phi (b_i b_{i+1} {b_i}^{-1}) \cdot \phi (b_i) \cdots \phi (b_m)$$
which is the same as $\phi (r_1) \cdots \phi (r_m)$.\\

\noindent From Propositions \ref{Project}, we are interested in the properties of the Hurwitz 
equivalence relation on factorizations in $S_n$. In $B_n$ we are interested in $\Delta ^2$ factorizations where all factors are powers of half-twists.\\
$\Delta ^2 = (\sigma _1...\sigma _{n-1})^n$ is a full $2 \pi $ twist of all the strings (in $B_n$)
and therefore, $\phi({\Delta ^2}) = 1_{S_n}$. If $H$ is a half twist by definition \ref{HT}, we get that $\phi (H) = (i,j) \quad 1 \geq i,j \geq n$. As a result, when projecting to $S_n$, we are interested in the properties of $1_{S_n}$ 
factorizations with transpositions or $1_{S_n}$ (when the power of the half twist is even) as factors.\\


\section{Hurwitz Equivalence Properties in $S_n$}


\bde
Let $\Gamma _1 \cdots \Gamma _m$, $\Gamma _i = (a_i, b_i) \quad 1 \leq a_i , b_i \leq n$ be a factorization. We define the graph of the factorization $G_F = (V_F , E_F ) $ where $V_F = \{ 1,...,n\} $ are the vertices and $E_F = \{(i,j) | \exists k \text{ s.t. } \Gamma _k =(i,j) \}$ are the edges of the factorization graph.
\ede


\bde
We define the weight of an edge $(i,j) \in E_F$ as the number of elements $\Gamma _k$ s.t. $\Gamma _k = (i,j)$. The weight of $(i,j)$ in the factorization $F$ will be noted as $W_F((i,j))$.
\ede


\noindent For a given graph $G_F$ we denote the graphs of its connected components as $G_F ^1,...,G_F ^d$ where $d$ is the number of connected components in the graph. For each connected component, let $G_F ^i = (V_F ^i, E _F ^i)$, where $V_F ^i$ are the vertices of $G_F ^i$ and $E_F ^i$ are the edges.\\


\bde
We define the weight of the connected component $G_F ^i$ to be:
$$\displaystyle W(G_F ^i) = \sum _{e_r\in E_F ^i} W(e_r)$$
\ede


\bthm \label{TheProposition}
Let $F_1, F_2$ be two $1_{S_n}$ factorizations with the same number of factors. Then 
$ F_1 \overset{HE}{\backsim} F_2$ if and only if $G_{F_1}$ and $G_{F_2}$ have the same number of connected components $G_{F_1} ^1,...,G_{F_1} ^d$ and $G_{F_2} ^1,...,G_{F_2} ^d$ respectively, and there exists a permutation $\pi$ s.t. $V_{F_1} ^i = V_{F_2} ^{\pi (i)}$ and $W(G_{F_1} ^i) = W(G_{F_2} ^{\pi (i)})$ for each $i \leq d$.
\ethm

\noindent In other words, two factorizations are Hurwitz equivalent if and only if the connected components of the factorizations graphs contain the same nodes and have the same weights.\\

\begin{figure}[h]
\begin{center}
\epsfxsize=8cm
\epsfysize=4cm
\epsfbox{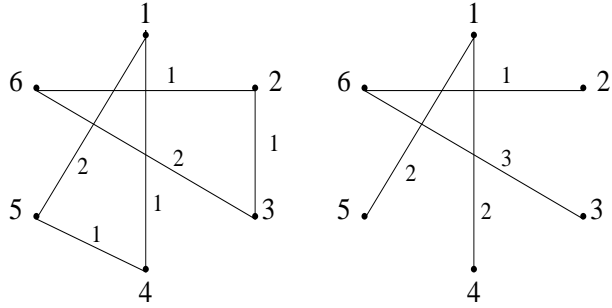}
\label{SnExample}
\caption{The graphs of the factorizations $F_1$ and $F_2$}
\end{center}
\end{figure}

\bex
\eexa
\noindent
The $1_{S_n}$ factorizations,\\
$$F_1 = (2,6) \cdot (1,4) \cdot (1,5) \cdot (3,6) \cdot (4,5) \cdot (1,5) \cdot (2,3) \cdot (3,6)$$ 
$$F_2 = (2,6) \cdot (1,5) \cdot (3,6) \cdot (3,6) \cdot (2,6) \cdot (1,5) \cdot (1,4) \cdot (1,4)$$ 
have connected components with the same nodes and and weights as shown in Figure \ref{SnExample}, and by Theorem \ref{TheProposition} they are Hurwitz equivalent.

\noindent The rest of the section will be devoted to the proof of Theorem \ref{TheProposition}. Starting with the first direction of the theorem.\\


\noindent
{\bf Proof of the first direction:} In the proof of the first direction we prove that if two factorizations are Hurwitz equivalent the factorizations have the same graph components with the same weights. Therefore, it is sufficient to show that when operating a single Hurwitz move, the vertices and weights of the graph's connected components will remain the same.\\

\noindent Let $F_1 = \Gamma _1 \cdots \Gamma _m$ and $F_2 = \Gamma _1 \cdots \Gamma _i \Gamma _{i+1} \Gamma _i ^{-1} \cdot \Gamma _i \cdots \Gamma _m$ the factorization obtained from $F_1$ by performing Hurwitz move $R_i$.\\
\noindent Let $\Gamma _j = (a_j, b_j), \quad 1 \leq a_j,b_j \leq n, \quad j\leq m$, then, in the cases where:
$$\{ a_i ,b_i \} \bigcap \{ a_{i+1}, b_{i+1} \} = \phi \text{    or}$$ 
$$\{ a_i, b_i \} \bigcap \{a_{i+1}, b_{i+1} \} = \{ a_i, b_i \}$$
\noindent We get that,
$$\Gamma _i \Gamma _{i+1} \Gamma _i ^{-1} = \Gamma _{i+1}$$
and the two factorizations have the same factors in a different order, so the factorizations graphs $G_{F_1}$ and $G_{F_2}$ are the same.\\

\noindent We are left with the case where $\Gamma _i =(a_i, b)$ and $\Gamma _{i+1}=(a_{i+1}, b)$.\\ 
$\Gamma _i \Gamma _{i+1} \Gamma _i ^{-1} = (a_i, a_{i+1})$ replaces $\Gamma _{i+1}$. 
We will show that the theorem still holds for this case.


\blem \label{FDLemma}
\elem
\begin{enumerate} 
\item If $v_1, v_2 \in V_{F_1} ^ r$ then $v_1, v_2 \in V_{F_2} ^ {r_1}$ for some $r_1$.
\item If $v_1 \in V_{F_1} ^ {r_1}$ and $v_2 \in V_{F_1} ^ {r_2}$, $r_1 \neq r_2$ then $v_1 \in V_{F_2} ^ {t_1}$ and $v_2 \in V_{F_2} ^ {t_2}$, for some $t_1,t_2$ s.t. $t_1 \neq t_2$.
\item If $V_{F_1}^ {r_1}= V_{F_2}^ {t_1}$ then  $W(G_{F_1}^ {r_1})= W(G_{F_2}^ {t_1})$.
\end{enumerate}


\noindent
{\bf Proof 1:} Since $v_1, v_2 \in V_{F_1} ^ r$ there is a sequence of edges connecting them in $ G_{F_1} ^r$, say, $\{ e_s \} _{s=1} ^p$.\\
\noindent
Since $G_{F_1} ^r$ is a connected component, $\{ e_s \} _{s=1} ^p \subset E_{F_1} ^ r$. So if $\Gamma _{i+1} \not\in E_{F_1} ^ r$ all  $\{ e_s \} _{s=1} ^p$ remain as elements in the factorization $F_2$, since only $\Gamma _{i+1}$ is replaced by $\Gamma _{i+1} \Gamma _i \Gamma _{i+1} ^{-1} $.\\
\noindent
In the case where  $\Gamma _{i+1} \in E_{F_1} ^ r$, every $s`$ s.t. $e_{s`} = \Gamma _{i+1}$ will be replaced in the sequence by $\Gamma _i \Gamma _{i+1} \Gamma _i ^{-1} = (a_i, a_{i+1})$ and $\Gamma _i = (a_i, b)$ which are elements in $F_2$ and connect $a_{i+1}$ with $b$.\\


\noindent 
{\bf Proof 2:} From (1), we conclude that the number of vertices in a connected component can only increase. 
Therefore, if $v_1 \in V_{F_1} ^ {r_1}$ and $v_2 \in V_{F_1} ^ {r_2}$, $r_1 \neq r_2$ and they are in the same 
connected component $V_{F_2} ^ {t}$ in $G_{F_2}$, then $ V_{F_1} ^ {r_1},  V_{F_1} ^ {r_2} \subset  V_{F_2} ^ {t}$. 
Therefore, there exists an edge $(v`_1, v`_2) \in E_{F_2} ^ {t}$ s.t. $v`_1, v`_2$ belongs to a different connected 
components in $G_{F_1}$. But the only edge that was added is $(a_{i+1}, a_i)$ and $a_{i+1}, a_i$ are in the same 
connected component in $G_{F_1}$.\\


\noindent 
{\bf Proof 3:} From (1) and (2), we see that the connected components remain the same. 
The weights of the connected components remain the same since all edges are the same except for $(a_i, b)$ 
that was replaced by $(a_i, a_{i+1})$. But $a_i, a_{i+1},b$ are all in the same connected component. 
Therefore the weight of all connected components remain the same.\\

\noindent
The Lemma proves that when performing a Hurwitz move on two transpositions the nodes of the connected components remain the same and so are the weights of the connected components. If one or both of the factors are $1_{S_n}$, the Hurwitz move does not change the factors only the order and therefore, the theorem still holds.\\

\noindent This concludes the proof of the first direction. \hfill $\qed $\\


\noindent
{\bf Proof of the second direction:} To complete Theorem \ref{TheProposition} we need to show that if two 
factorizations have the same connected components with the same weights they are Hurwitz equivalent. To prove that 
we will show that each factorization is Hurwitz equivalent to a standard canonical factorization which depends only on 
the nodes of the factorization`s connected components and their weights.\\


\blem \label{ConjPerm}
Let $(a,b) \cdot (c,d)$ be a factorization is $S_n$ then,
\begin{enumerate}
\item By performing Hurwitz move $R_0$ we get:\\
$$(a,b) \cdot (c,d) \overset{HE}{\backsim}
{	\begin{cases} 
		{(c,d) \cdot (a,b),  \text{ if } \{ a,b\} \bigcap \{ c,d \} = \phi  }\\
		{(c,d) \cdot (a,b),  \text{ if } \{ a,b\} \bigcap \{ c,d \} = \{ a,b\}  }\\
		{(a,d) \cdot (a,b),  \text{ if } b = c \text{ and } a \neq d  }
\end{cases} } $$
\item By performing Hurwitz move ${R_0}^{-1}$ we get:\\
$$(a,b) \cdot (c,d) \overset{HE}{\backsim}
{	\begin{cases} 
		{(c,d) \cdot (a,b),  \text{ if } \{ a,b\} \bigcap \{ c,d \} = \phi  }\\
		{(c,d) \cdot (a,b),  \text{ if } \{ a,b\} \bigcap \{ c,d \} = \{ a,b\}  }\\
		{(c,d) \cdot (a,d),  \text{ if } b = c \text{ and } a \neq d  }
\end{cases} } $$
\item $(a,b) \cdot 1_{S_n} \overset{HE}{\backsim} 1_{S_n} \cdot (a,b)$
\end{enumerate}
\elem
\noindent
{\bf Proof:} Trivial.\\

\noindent
From Lemma \ref{ConjPerm} (3) the $1_{S_n}$ factors commutes with all other factors. Each factorization is Hurwitz equivalent to a factorization where all $1_{S_n}$ factors are on the left of the factorization and the two factorizations have the same transpositions as factors. Since the theorem requires that the weights of the connected components is equal and that both factorizations have the same number of factors, the number of $1_{S_n}$ factor is equal.\\
\noindent Therefore, to prove Theorem \ref{TheProposition} we can ignore the $1_{S_n}$ factors, and find standard canonical form to the transposition factors only.\\

\noindent Let $F_1 = \Gamma _1\cdots \Gamma _m$ be a factorization of $1_{S_n}$ where all factors are of transpositions.


\blem \label{Commutes}
If $\Gamma _j \not\in E _{F_1} ^{t_1}$ and $\Gamma _{j+1} \not\in E _{F_1} ^{t_2}$, $t_1 \neq t_2$ then 
$\Gamma _1 \cdots \Gamma _j \cdot \Gamma _{j+1} \cdots \Gamma _m \overset{HE}{\backsim} \Gamma _1 \cdots \Gamma _{j+1} \cdot \Gamma _{j} \cdots \Gamma _m$.
\elem
\noindent 
{\bf Proof:} Since $\Gamma _j$, $\Gamma _{j+1}$ belong to a different connected component, they do not connect the same 
vertex and therefore, $\Gamma _j \Gamma _{j+1} \Gamma _j ^{-1}= \Gamma _{j+1}$ (See Lemma \ref{ConjPerm}). Therefore, by operating 
Hurwitz move $R_i$ they commute.\\

\noindent As a result, for each connected component, all elements of the component commutes with all elements of other 
components. Therefore, factorization is Hurwitz equivalent to a factorization with the same factors ordered according 
to the component they belong too. For example, order the connected components by the lowest vertex they 
contain, then gather all factors of the first component to the left, and after them the factors of the second component and so on.\\

\noindent Let $\{ G _{F_1} ^r \} _{r=1} ^s$ be the distinct connected components of $G_{F_1}$. 
From Lemma \ref{Commutes}, $F_1 \overset{HE}{\backsim} f_1 \cdots f_s$ where $f_r$ is a factorization with elements from $E _{F_1} ^r$. 
The length of the factorization $f_r$ is equal to $W(G_{F_1} ^r)$ and $s$ is the number of connected components.
Therefore, to conclude the proof, it is sufficient to show that each 
$f_r$ is Hurwitz equivalent to a standard canonical factorization which depends only on the length of $f_r$ (which can never be 
changed by Hurwitz moves) and $V _{F_1} ^r$.

\noindent We define an order on $V _{F_1} ^r$ vertices, $V _{F_1} ^r = \{ v_{t_1},...,v_{t_l} \} $. Note that since 
$F_1 = 1_{S_n}$ then, $f_r = 1_{S_n}$ as a product.\\


\blem \label{LemmaLeft}
Let $f=\Gamma _1 \cdots \Gamma _m $ be a factorization with a single connected component, $G_f ^1$, then 
$\forall v_1,v_2 \in V_f ^1, \quad f \overset{HE}{\backsim} (v_1,v_2) \cdot \gamma _1 \cdots \gamma _{m-1}$.
\elem

\noindent
{\bf Proof:} Proof by induction on the minimal length of the path connecting $v_1$ with $v_2$. In the case where 
the minimal length is $1$, there exists $1 \leq j \leq m$ s.t. $\Gamma _j = (v_1, v_2)$. Operating $\{ {R_k}^{-1}\} _{k=j-2} ^0$ 
sequence of Hurwitz moves, we get a factorization $(v_1,v_2) \cdot \gamma _1 \cdots \gamma _{m-1}$ 
(See Lemma \ref{ConjPerm}). We will assume that the lemma is true for a path with length less than $n$, we will prove 
that the factorization where the minimal path between $v_1$ and $v_2$ is $n$, is Hurwitz equivalent to a factorization 
which the path between $v_1$ and $v_2$ is of length $n-1$:\\

\noindent Let $(a_1,a_2),(a_2, a_3),...,(a_n, a_{n+1})$ be the minimal path between $v_1 = a_1$ and $v_2 = a_{n+1}$. To prove the 
above we will perform another induction, on the number of factors which are in between $(a_1, a_2)$ and $(a_2, a_3)$. We 
will assume that $(a_1, a_2)$ is left to $(a_2, a_3)$: 
$$f= \cdots (a_1,a_2) \cdot (b_1, c_1) \cdot (b_2, c_2) \cdots (b_k, c_k)\cdot (a_2, a_3)\cdots $$
Let $k$ be the number of factors between $(a_1, a_2)$ and $(a_2, a_3)$.\\
In the case where $k=0$, $(a_1,a_2) \cdot (a_2, a_3) \overset{HE}{\backsim} (a_1,a_3),(a_1, a_2)$ and we are done since the new 
factorization contains the path $(a_1,a_3),(a_3, a_4),...,(a_n, a_{n+1})$ which is of length $n-1$.\\
In the case where $k > 0$:

\begin{emphit} 
\item If $\{ a_1, a_2 \} \bigcap \{ b_1, c_1 \} = \phi $ then, 
$(a_1,a_2) \cdot (b_1, c_1) \overset{HE}{\backsim} (b_1,c_1) \cdot (a_1, a_2)$ (By Lemma \ref{ConjPerm}), and now $(a_1,a_2)$ 
and $(a_2, a_3)$ are separated by $k-1$ factors, and so we are done.

\item If $a_1 =  b_1$ then, 
$(a_1,a_2) \cdot (a_1, c_1) \overset{HE}{\backsim} (a_2,c_1) \cdot (a_1, a_2)$ (By Lemma \ref{ConjPerm}), and now $(a_1,a_2)$ 
and $(a_2, a_3)$ are separated by $k-1$ factors. Note that $(a_1, c_1)$ is not an element in the path 
(since the path is minimal).

\item If $a_2 =  b_1$ then, 
$(a_1,a_2) \cdot (a_2, c_1) \overset{HE}{\backsim} (a_1,c_1) \cdot (a_1, a_2)$ (By Lemma \ref{ConjPerm}), and now $(a_1,a_2)$ 
and $(a_2, a_3)$ are separated by $k-1$ factors. Note that if $(a_2, c_1)$ is in the path then $c_1 = a_3$ and then 
$k=0$.
\end{emphit}
\noindent This concludes the proof of Lemma \ref{LemmaLeft}.


\blem \label{DoubleMoves} 
\elem
\begin{enumerate}
\item $(a,b)\cdot (a,b)\cdot (a,c) \cdot (a,c) \overset{HE}{\backsim} (a,c) \cdot (a,c) \cdot (a,b)\cdot (a,b)$.
\item $(a,b)\cdot (a,b)\cdot (a,c) \cdot (a,c) \overset{HE}{\backsim} (a,b) \cdot (a,b) \cdot (b,c)\cdot (b,c)$.
\end{enumerate}

\noindent
{\bf Proof 1:} $(a,b)\cdot (a,b)\cdot (a,c) \overset{HE}{\backsim} (a,c) \cdot (a,b)\cdot (a,b)$ By operating Hurwitz moves $R_1$ and $R_0$ and 
therefore, $(a,b)\cdot (a,b)\cdot (a,c) \cdot (a,c) \overset{HE}{\backsim} (a,c) \cdot (a,c) \cdot (a,b)\cdot (a,b)$.\\\\
{\bf Proof 2:} By performing the Hurwitz moves $R_1 ^{-1}, R_2 ^{-1}, R_1 ^{-1}$.\\

\noindent Now we are ready to start forming $f_r$ into a standard canonical form:\\
\noindent $v_{t_1}, v_{t_2} \in V_f ^r$, from Lemma \ref{LemmaLeft}, 
$$f_r \overset{HE}{\backsim} (v_{t_1}, v_{t_2}) \cdot f_r ^1$$ 
where $f_r ^1$ is the factorization with the $W(G_{f_r} ^1)-1$ other factors.\\
\noindent 
$f_r ^ 1 = (v_{t_1}, v_{t_2})$ since $f_r = 1_{S_n}$ and $f_r ^ 1 = (v_{t_1}, v_{t_2})^{-1} f_r$.\\
\noindent
Because $f _r ^1 = (v_{t_1}, v_{t_2})$, $f_r ^ 1$ contains a path connecting $v_{t_1}$ with $v_{t_2}$. Again, using Lemma \ref{LemmaLeft} we get, 
$$f_r \overset{HE}{\backsim} (v_{t_1}, v_{t_2})\cdot (v_{t_1}, v_{t_2}) \cdot f_r ^2 \text{ and } f_r ^2 = 1_{S_n}$$\\

\noindent By the first direction of Theorem \ref{TheProposition} 
$(v_{t_1}, v_{t_2})\cdot (v_{t_1}, v_{t_2}) \cdot f_r ^2$ still
creates a single connected component. Therefore, there is a path between $v_{t_2}$ and $v_{t_3}$, 
which means that $f_r ^2$ contains a path from $v_{t_2}$ to $v_{t_3}$ or from $v_{t_1}$ to $v_{t_3}$ 
(for example in some cases where the path in $G_f ^r$ includes $ (v_{t_1},v_{t_2})$).\\
In the first case we get 
$f_r \overset{HE}{\backsim} (v_{t_1}, v_{t_2})(v_{t_1}, v_{t_2})(v_{t_2}, v_{t_3})(v_{t_2}, v_{t_3}) \cdot f_r ^4$ 
and in the second case we get
$f_r \overset{HE}{\backsim} (v_{t_1}, v_{t_2})(v_{t_1}, v_{t_2})(v_{t_1}, v_{t_3})(v_{t_1}, v_{t_3}) \cdot f_r ^4$ 
which is by Lemma \ref{DoubleMoves} Hurwitz equivalent to the first case.\\

\noindent We continue with this process to bring $f_r$ to the form:\\ 
$(v_{t_1}, v_{t_2})\cdot (v_{t_1}, v_{t_2})\cdot (v_{t_2}, v_{t_3}) \cdot (v_{t_2}, v_{t_3}) \cdots (v_{t_{m-1}}, v_{t_m})(v_{t_{m-1}}, v_{t_m}) \cdot f_r ^{2m-2}$\\
\noindent Assume we came to the point where:\\
$f_r \overset{HE}{\backsim} (v_{t_1}, v_{t_2}) \cdot (v_{t_1}, v_{t_2}) \cdots (v_{t_{k-1}}, v_{t_k})\cdot (v_{t_{k-1}}, v_{t_k}) \cdot f_r ^{2k-2}$ 
where $k < m$ and $f_r ^{2k-2}$ is a factorization with $W(G_f ^r)-2k +2$ factors.\\
\noindent Same as before, the new factorization creates a connected graph and $f_r ^{2k-2} = 1_{S_n}$. 
Since the graph is connected, $f_r ^{2k -2}$ contains a path from $v_{t_{k+1}}$ to one of the vertices $v_{t_s} \quad s \leq k$. 
So, there is a path from $v_{t_s} \quad (s \leq k)$ to $v_{t_{k+1}}$ which does not include the 
factors left to $f_r ^{2k-2}$ because they create a connected graph which does not include $v_{t_{k+1}}$.
So, from Lemma \ref{LemmaLeft}, 
$$f_r \overset{HE}{\backsim} (v_{t_1}, v_{t_2}) \cdot (v_{t_1}, v_{t_2}) \cdots (v_{t_{k-1}}, v_{t_k})\cdot (v_{t_{k-1}}, v_{t_k}) (v_{t_{k+1}}, v_{t_s})\cdot f_r ^{2k -1} \text{,}$$ 
and since $f_r ^{2k -1} = (v_{t_{k+1}}, v_{t_s})$ as a product, there is a path from $v_{t_{k+1}}$ to $v_{t_s}$.\\
\noindent By Lemma \ref{LemmaLeft}, 
$$f_r \overset{HE}{\backsim} (v_{t_1}, v_{t_2}) \cdot (v_{t_1}, v_{t_2}) \cdots (v_{t_{k-1}}, v_{t_k})\cdot (v_{t_{k-1}}, v_{t_k}) \cdot (v_{t_{k+1}}, v_{t_s})\cdot (v_{t_{k+1}}, v_{t_s})\cdot f_r ^{2k} \text{.}$$\\
Now, only using the factors left to $f_r ^{2k}$ we need to change $(v_{t_{k+1}}, v_{t_s})\cdot (v_{t_{k+1}}, v_{t_s})$ 
to $(v_{t_{k}}, v_{t_{k+1}})\cdot (v_{t_{k}}, v_{t_{k+1}})$. Since $s\leq k$ there is a path from $v_{t_s}$ to $v_{t_k}$  
in the graph created by the factors on the left, from this fact and using Lemma \ref{DoubleMoves} we see that the factorizations:\\
$(v_{t_1}, v_{t_2}) \cdot (v_{t_1}, v_{t_2}) \cdots (v_{t_{k-1}}, v_{t_k})\cdot (v_{t_{k-1}}, v_{t_k}) \cdot (v_{t_s}, v_{t_{k+1}})\cdot (v_{t_s}, v_{t_{k+1}}) \text{ and}$\\
\noindent
$(v_{t_1}, v_{t_2}) \cdot (v_{t_1}, v_{t_2}) \cdots (v_{t_{k-1}}, v_{t_k})\cdot (v_{t_{k-1}}, v_{t_k}) \cdot (v_{t_k}, v_{t_{k+1}})\cdot (v_{t_k}, v_{t_{k+1}})$\\
are Hurwitz equivalent, since Lemma \ref{LemmaLeft} allows us to commute couples of transpositions, or to change one vertex in the couple if the two couples have a common vertex.

\noindent To complete the proof of Theorem \ref{TheProposition} we need to show that we can also bring the right factors, 
$f_r ^{2m-2}$ to a standard form. This can be done in a similar way to the above procedure:\\
\noindent Take the first factor in $f_r ^{2m-2}$, i.e. $f_r ^{2m-2} = (v_{t_x}, v_{t_y}) \cdot f_r ^{2m-1}$ and again by 
using Lemma \ref{LemmaLeft} we get $f_r ^{2m-2} \overset{HE}{\backsim} (v_{t_x}, v_{t_y}) \cdot (v_{t_x}, v_{t_y}) f_r ^{2m}$.\\

\noindent Using Lemma \ref{DoubleMoves} and the factors on the left, we can change $(v_{t_x}, v_{t_y}) \cdot (v_{t_x}, v_{t_y})$ 
to $(v_{t_1}, v_{t_2}) \cdot (v_{t_1}, v_{t_2})$ since the graph of the factors on the left is the same as $f_r$.
This concludes the proof of Theorem \ref{TheProposition} since every factorization $f_r$ is Hurwitz equivalent to:\\
$$(v_{t_1}, v_{t_2}) \cdot (v_{t_1}, v_{t_2}) \cdots (v_{t_{m-1}}, v_{t_m})\cdot (v_{t_{m-1}}, v_{t_m}) (v_{t_1}, v_{t_2}) \cdots (v_{t_1}, v_{t_2})$$ 
Which depends only on the factorization graph and the number of factors in the factorization. \hfill $\qed $

\end{document}